\documentclass[12pt]{amsart}

\usepackage{hyperref, color}

\newtheorem{theorem}{Theorem}[section]
\newtheorem{lemma}[theorem]{Lemma}

\newtheorem{corollary}[theorem]{Corollary}

\theoremstyle{remark}
\newtheorem{remark}[theorem]{Remark}

\numberwithin{equation}{section}

\newcommand{\z}{\langle}
\newcommand{\x}{\rangle}

\DeclareMathOperator{\indw}{Ind_w}

\DeclareMathOperator{\diver}{div}
 
\begin{document}

\title[Index Estimates for Free Boundary CMC Surfaces]
{Index Estimates for Free Boundary Constant Mean Curvature Surfaces}
\author[Cavalcante]{Marcos P. Cavalcante}
\address{ \newline 
Instituto de Matem\'atica 
\newline Universidade Federal de Alagoas
\newline Campus A. C. Sim\~oes, BR 104 - Norte, Km 97, 57072-970.
\newline Macei\'o - AL -Brazil}
\email{marcos@pos.mat.ufal.br}
\thanks{The first author was partially  supported  by CNPq-Brazil.}

\author[de Oliveira]{Darlan F. de  Oliveira}
\address{
\newline Departamento de Ci\^encias Exatas 
\newline  Universidade Estadual de Feira de Santana 
\newline Avenida Transnordestina, S/N, Novo Horizonte, 44036-900
 \newline Feira de Santana - BA - Brazil}
\email{darlanfdeoliveira@gmail.com}


\date{\today}


\begin{abstract} {
In this paper we consider compact constant mean curvature surfaces with boundary immersed in a 
mean convex body of the Euclidean  space or in the unit sphere. 
We  prove that the Morse 
index is bounded from below by a linear function of the genus and number of boundary components.}
\end{abstract}

\maketitle

\section{Introduction}

Let $W$ be a Riemannian manifold with nonempty boundary
such that its  boundary $\partial W$ is a union of smooth hypersurfaces. 
Let  $M\subset W$ be a  compact  constant mean curvature hypersurface 
such that $M$ intersects the regular part of $\partial W$ 
along its boundary in a right angle. 
It is well known that such hypersurfaces are critical points of the area functional 
for variations of $M$ that preserve the enclosed volume and keep the boundary freely 
on $\partial W$.
These hypersurfaces arise in many geometrical and physical problems and 
are referred as \emph{free boundary CMC hypersurfaces}
(FBCMC hypersurfaces, for short). They have been
studied since the 19th century and  still now form  a very active topic in differential geometry. 
We refer the  reader to the books of Finn \cite{Finn} and L\'opez \cite{Lopez}, as well as the references therein, for a nice introduction to this subject.  

An important problem about FBCMC hypersurfaces is to classify those ones that are 
\emph{stable},
that is, whose the second variation of the area is nonnegative for volume preserving variations. 
For instance, in the case that $W$ is a geodesic ball in a space form, 
a well-known conjecture asserted that spherical caps are the only solutions. 
It was confirmed by the works of Ros-Vergasta  \cite{RV} and Nunes \cite{Nunes} in dimension two and more recently by Wang and Xia \cite{WangXia} in any dimension.
Other results on stable FBCMC hypersurfaces can be found for instance in 
\cite{AinouzSouam}, \cite{Athanassenas}, \cite{Barbosa}, \cite{ChoeKoiso}, \cite{LiXiong2}, \cite{LiXiong}, 
\cite{Lopez14}, \cite{R08} and \cite{RosSouam}.

When $M$ is a stable FBCMC surface immersed in a mean convex region 
$W\subset\mathbb R^3$  Ros showed in  \cite[Theorem 9]{R08} that there are just a few 
possibilities for the genus and the number of boundary components of $M$
(see Corollary \ref{c1} below).

When $M$ is not stable, the Jacobi operator $J$ associated to the second variation of the area has non zero \emph{Morse index} (see Section \ref{preliminaries} for precise definitions).  
Geometrically, the index of $M$ is the number of directions whose admissible 
variations decrease area. It will be denoted by $\indw(M)$.

The  case of free boundary minimal surfaces is of special interest and  the 
works of Fraser and Schoen \cite{FS11, FS13, FS16} have motivated many research in this case.
For free boundary minimal hypersurfaces all variations keeping the boundary freely in the boundary are allowed, not only volume preserving variations.   
In this setting, recently  Ambrozio, Carlotto and Sharp in \cite{ACS} and  Sargent in \cite{PS},
proved independently that if $M$ is a free boundary minimal surface immersed in a
mean convex region $W\subset \mathbb R^3$ then the index is bounded from below by
$(2g + k-1)/3$, where $g$ is the genus of $M$ and $k$ is the number of boundary components.
In higher dimensions, they  also obtained  lower bounds for the index in terms of the
dimension of the  first relative homology group with real coefficients. 
The technique presented in these results uses the coordinates of harmonic forms as test functions
and is inspired by
previous works on eigenvalue estimates  and index estimates for 
minimal hypersurfaces without boundary 
(see \cite{R06}, \cite{Savo} and  \cite{ACS}). 

Following these lines, in this paper, we obtain a comparison theorem for the eigenvalues
of $J$ in terms of the mean curvature and the eigenvalues of the Hodge Laplacian $\Delta$ 
acting on 1-forms and, as a by-product, we obtain lower bounds for the Morse index
of FBCMC
surfaces in mean convex regions of the Euclidean space or the unit sphere. More precisely,
our results are the following.

\begin{theorem}\label{esp}
Let $W$ be a region  of $\mathbb{R}^3$  such that its boundary is a union of 
smooth mean convex surfaces, and let $M^2$ be a compact, 
orientable, FBCMC surface immersed in the mean convex side of
$W$ and whose boundary intersects the regular part of $\partial W$. 
Then for all positive integers $\alpha$ we have
\begin{eqnarray*} \lambda_{\alpha}^J\leq -2H^2+\lambda_{m(\alpha)}^{\Delta},
\end{eqnarray*}
where $m(\alpha)>6(\alpha - 1).$

In particular, if $M$ has genus $g$ and $k$ boundary components. Then,
\[
\indw(M)\geq\frac{2g+k-4}{6}.
\]
\end{theorem}

As an immediate consequence we obtain the result of Ros \cite{R08} cited above.
\begin{corollary}[Ros]\label{c1} Under the conditions of Theorem \ref{esp}, if $M$ is stable, then only possibilities 
for $g$ and $k$ are 
\begin{enumerate}
\item $g=0 $ and $k\leq 4$;
\item $g=1$ and $k=1$ or $2$.
\end{enumerate}
\end{corollary}

In the case of FBCMC surfaces immersed in mean convex domains of $\mathbb S^3$, the result read as follows. 

\begin{theorem}\label{espS3}
Let $W$ be a region  in the unit sphere $\mathbb{S}^3$  such that its boundary is a union of 
smooth mean convex surfaces, and let $M^2$ be a compact, 
orientable, FBCMC surface immersed in the mean convex side of 
$W$ and whose boundary intersects the
regular part of $\partial W$. 
Then for all positive integers $\alpha$ we have
\begin{eqnarray*} \lambda_{\alpha}^J\leq -2(H^2+1)+\lambda_{m(\alpha)}^{\Delta},
\end{eqnarray*}
where $m(\alpha)>8(\alpha-1).$

In particular,
\[
\indw(M)\geq\frac{2g+k-5}{8}.
\]

\end{theorem}

\begin{corollary} Under the conditions of Theorem \ref{espS3}, 
if $M$ is stable, then only possibilities 
for $g$ and $k$ are 
\begin{enumerate}
\item $g=0 $ and $k\leq 4$;
\item $g=1$ and $k=1$ or $2$;
\item $g=2$ and $k=1$.
\end{enumerate}
\end{corollary}

This paper is organized as follows. In Section \ref{preliminaries} we present somen definitions
and basic results to be used in the proofs.  Section \ref{test} is devoted to computing the
Jacobi operator of the test functions given by the coordinates of harmonic forms. 
In Section \ref{proofs} we present the proof of Theorem \ref{esp}. The proof of Theorem
\ref{espS3} is analogous and will be omitted.

\section{Preliminaries}\label{preliminaries}

Let us denote by $W$ a domain of the Euclidean space  $\mathbb R^3$ not
necessarily compact. For simplicity, let us assume that $W$ has smooth boundary and
let $\nu$ be the unit normal vector field along $\partial W$.
We recall that the second fundamental form and the mean curvature of $\partial W$ 
with respect to $\nu$ are defined respectively by
\[
II^{\partial W}(X,Y) = \langle -D_X\nu, Y\rangle, \quad \textrm{ for }X,Y\in T\partial W,
\]
and
\[
H^{\partial W} = \frac 1 2 \textrm{tr}II^{\partial W},
\]
where $D$ is the Levi-Cevita connection in the Euclidean space. 
If we can choose the unit normal  $\nu$ along $\partial W$ 
such that $\nu$ points outward $W$ and $II^{\partial W}$ is negative defined, then
we say that $W$ is \emph{convex}. 
If it can be chosen such that  $H^{\partial W}\leq 0$, then we say that
$W$ is \emph{mean convex}. 

Let $x:M\to W$ be a compact oriented surfaced with boundary which is  
properly immersed, that is,  $M \cap  W = \partial M$.
Fixed a unit normal vector field $N$ along $x$ we denote by $A$  the shape operator 
associated to the second fundamental form $II^{\partial M}$
of $M$ with respect to $N$, namelly 
\[AX=-D_XN, \quad \textrm{ for }X\in TM,
\] 
and
\[
II^{\partial M}(X,Y) = \langle AX, Y\rangle, \quad \textrm{ for }X,Y\in TM,
\]
We say that $M$ is free boundary if $\partial M$  meets $\partial W$ orthogonally.
 
From now on, let us assume that  $W$ is a mean convex
domain of $\mathbb R^3$ and $M$ is a 
\emph{free boundary constant mean curvature surface} properly immersed in $W$.
Such surfaces are critical points of the area functional  for variations that
preserves the enclosed volume (see \cite{N85}). 
It is easy to see that, up a parametrization, we can consider only normal variations,
say generated by $X=uN$, where $u\in \mathcal F$ and
\[
\mathcal F = \{u:M\to \mathbb R: u \textrm{ is smooth up to the boundary 
and}\, \int_M u\, dM=0\}.
\]
In this setting, the second variation of area functional is given by the following 
quadratic form (see \cite{RV},  \cite{R08})
\[Q(u,u)=\int_M(u\Delta u - \|A\|^2u^2)\, dM+
\int_{\partial M}(u \eta(u)+II^{\partial W}(N,N)u^2)\,ds.
\]
Here $\eta$ is the unit conormal vector to $M$  along its boundary.
Note that the free boundary condition means that $\eta=\nu$ along $\partial W$.
We point out that in this paper we are using the geometric definition of the Laplacian
operator, that is, $\Delta u = \diver\nabla u$, where $\diver X = -\textrm{tr} \nabla X$.

The \emph{index} of $M$, denoted by $\indw(M)$, is 
defined as the maximal dimension of a 
subspace of $\mathcal{F}$ on which $Q$ is negative. 
Geometrically, the index indicates the number of directions whose variations decrease area. 
In particular, we say that $M$ is \emph{stable} if the index is zero. 
We also can define the index as the number of negative eigenvalues of the boundary problem
associated to the quadratic form $Q$. More precisely, 
let us consider $J=\Delta-\|A\|^2$ the
Jacobi operator of $M$. 
We say that $\lambda^J\in \mathbb R$ is an eigenvalue of $Q$ if there exists a nontrivial 
eigenfunction $u\in \mathcal{F}$ such that
\begin{eqnarray}\label{proauto}
\left\{\begin{array}{l}
Ju=\lambda^J u, \quad \mbox{ in} \ M,\\
\frac{\partial u}{\partial \eta}=-II^{\partial B}(N,N)u, \quad \mbox{on} \ \partial M. 
\end{array}\right.
\end{eqnarray}

It is well known that there exists a non-decreasing sequence 
\[
\lambda_1<\lambda_2\leq\cdots\leq \lambda_k\leq\cdots \nearrow \infty
\]
of eigenvalues  associated to a $L^2(M)$-orthonormal basis 
$\{\phi_1,\cdots,\phi_k,\cdots\}$ of solutions to the eigenvalue problem (\ref{proauto}) and  satisfying the min-max characterization 
\[
\lambda_k^J=\min_{u\in \mathcal{G}}\frac{Q(u,u)}{\int_Mu^2dM},
\]
where $\mathcal{G}=\langle\phi_1,\cdots,\phi_{k-1}\rangle^{\bot}\setminus \{0\}.$
It follows immediately that $\indw(M)$ equals the number of negative eigenvalues of the
problem (\ref{proauto}).

 In order to give lower bounds for the index of $M$ in terms of its topological invariants
we will construct admissible eigenfunctions in $\mathcal F$ using harmonic vector fields,
or equivalently harmonic $1$-forms.
Let us consider the set of $1$-forms that are normal at $\partial M,$ that is, 
\[
\Omega^1(M,\partial M):=\{w\in\Omega^1(M),i^*w=0\},
\]
where $i:\partial M\to M$ is the inclusion map.

We then consider the space of tangential harmonic 
$1$-forms 
\[
\mathcal{H}^1_N(M):=\{w\in\Omega^1(M,\partial M);dw=0 \ \mbox{and} \ \delta w=0\}
\]
and  the space of normal harmonic  $1$-forms 
\[
\mathcal{H}^1_T(M):=\{\star w\in\Omega^1(M,\partial M);dw=0 \ \mbox{and} \ \delta w=0\}.
\] 

Above  $d$ is the  exterior derivative operator and 
$\delta $ is the  interior derivative operator defined by 
$\delta=-\star d\star$, where 
$\star:\Omega^1(M)\to\Omega^{1}(M)$ is the Hodge star operator.

It is well known that these spaces are closed related to the topology of the underline 
manifold. In fact we have the following result (see \cite{ACS} or \cite{PS}).
\begin{lemma}\label{dim}
Let $M^2$ be a compact, orientable surface with non-empty boundary $\partial M.$ If $M$ has genus $g$ and $k\geq1$ boundary components, then
\[
\dim \mathcal{H}^1_T(M)=\dim  \mathcal{H}^1_N(M)=2g+k-1.
\]
\end{lemma}
%
\section{Test functions and harmonic vector fields} \label{test}


Denoting by   $\mathcal{E}=\{\bar E_1,\bar E_2,\bar E_{3}\}$  
the canonical basis  in $\mathbb{R}^{3}$ we will considere
${E_i}:=\bar E_i-\langle \bar E_i,N\rangle N$, 
the orthogonal projection of $\bar E_i$ on $TM$. 
We also consider the smooth support functions 
$g_i:M\to\mathbb R$,  $g_i:=\left\langle \bar{E_i},N\right\rangle,$
for $1\leq i\leq 3.$

Given a smooth vector field  $\xi\in TM$  on $M$ we will use its coordinates
and the coordinates of $\star \xi$
as test functions. Namely, for each $1\leq i\leq 3,$  we define
$w_i,\, \bar w_i:M\to \mathbb R$ as
\[
w_i:=\langle  E_i,\xi\rangle \textrm{ and }\,
\bar w_i:=\langle  E_i,\star \xi\rangle.
\]

In order to compute the Jacobi operator of $w_i$ and $\bar w_i$ we recall the following lemma of
local nature proved in \cite{CdO} (see also \cite{R07}). 
\begin{lemma}\label{lapwi}Let $M^{2}$ be an orientable CMC  surface in $\mathbb{R}^3$. 
Then, using the above notation we have 
\[
\Delta w_i=(\|A\|^2-4H^2)w_i+2H\langle AE_i,\xi\rangle-2g_i\langle A,\nabla \xi\rangle+\langle E_i,\Delta \xi\rangle,
\]
and
\[
\Delta \bar w_i=(\|A\|^2-4H^2)\bar w_i+2H\langle AE_i,\star \xi\rangle-
2g_i\langle A,\nabla \star\xi\rangle+\langle E_i,\Delta \star\xi\rangle,
\]
for $1\leq i\leq 3$. 
\end{lemma}

Now we note  that when the vector field $\xi$  is harmonic and tangential along $\partial M$
its coordinates are  admissible functions to compute the index of CMC surfaces. 
More precisely we have:

\begin{lemma}\label{wi} If $\xi\in TM$ is a harmonic vector field which is tangential 
in $\partial M,$ then $w_i\in\mathcal F$, that is, 
\[
\int_M w_i \, dM=0, 
\]
for $1\leq i\leq 3.$
\end{lemma}
\proof Note that $E_i=\nabla x_i$, $1\leq i\leq 3,$ where 
 $x=(x_1,x_2,x_3): M\to B$ is the immersion map.
 Then we have 
\begin{eqnarray*}		
\int_M w_i\, dM&=&\int_M\langle \nabla x_i,\xi\rangle\, dM\\
                 &=&\int_M x_i\diver\xi\, dM+\int_{\partial M}x_i\langle \xi,\eta\rangle \,ds=0.
\end{eqnarray*}
In fact, 
$\diver\xi=0$ since $\xi$ is harmonic   and $\langle\xi,\eta\rangle=0$, 
since $\xi$ tangential to $\partial M.$
\endproof
\begin{remark}
In general the functions  $\bar w_i$, $1\leq i\leq 3$, have not
mean value zero.
However, we will see in Section \ref{proofs} 
that  if $\dim \mathcal{H}^1_T(M)$ is  large enough then we can choose $\xi$ 
such that 
\[
\int_M \bar w_i \, dM=0, 
\]
for $1\leq i\leq 3.$
\end{remark}

We conclude this section computing the boundary term of the quadratic form $Q$
on $w_i$ and  $\bar w_i,$. 
\begin{lemma}\label{wi2} 
If $\xi\in TM$ is a harmonic vector field which 
is tangential in $\partial M,$ then 
\begin{equation}\label{bordo1}
\sum_i\int_{\partial M}(w_i\eta(w_i)+II^{\partial W}(N,N)w_i^2)\, ds
=2\int_{\partial M}H^{\partial W}\| \xi\|^2\, ds
\end{equation}
and
\begin{equation}\label{bordo2}
 \sum_i\int_{\partial M}(\bar w_i\eta(\bar w_i)+II^{\partial W}(N,N)\bar w_i^2)\, ds
=2\int_{\partial M}H^{\partial W}\| \xi\|^2\, ds.
\end{equation}
\end{lemma}
\proof 

We first note that for any vector field $X$ we have 
\begin{eqnarray*}
\langle \nabla _\eta E_i, X\rangle  &=&
\eta\langle \bar E_i,X\rangle -\langle\bar E_i,\nabla _\eta X \rangle\\
&=&\langle \bar E_i, D_\eta X - \nabla _\eta X\rangle\\
&=& \langle \bar E_i,N\rangle\z \xi,A\eta\x.
\end{eqnarray*}
For $\xi$ tangential at $\partial M$ we have 
{\setlength\arraycolsep{1pt}
\begin{eqnarray*}
\sum_i\int_{\partial M}w_i\eta(w_i)\,ds
&=&\sum_i\int_{\partial M}w_i(\z\nabla_{\eta}\xi,E_i\x+\langle \bar E_i,N\rangle\z \xi,A\eta\x)\,ds\\
&=&\int_{\partial M}\z\nabla_{\eta}\xi,\xi\x\,ds\\
&=&\int_{\partial M}\z\nabla_{\xi}\xi,\eta\x\,ds\\
&=&-\int_{\partial M}\z\nabla_{\xi}\eta,\xi\x\,ds\\
&=&\int_{\partial M}II^{\partial W}(\xi,\xi)\,ds.
\end{eqnarray*}}
Since $\xi$ and $N$ form an orthogonal basis of the tangent space 
of $\partial W$ along $\partial M$ we conclude the proof by noting 
that  $$II^{\partial W}(\xi,\xi)+II^{\partial W}(N,N)\|\xi\|^2=2H^{\partial W}\|\xi\|^2 .$$
The proof of assertion (\ref{bordo2}) follows the same steps above noting additionally that 
the Levi-Civita connection $\nabla$ commutes with the Hodge star operator 
$\star$.
\endproof
%

\section{Proof of Theorem \ref{esp}}\label{proofs}

The proofs follow the same spirit as we did in \cite{CdO} but taking in account the 
boundary term. We point out here that we can assume without loss of generality that
$\partial W$ is smooth.  
Let $\xi_1, \xi_2, \ldots, \xi_m$ be 
the first $m$ eigenfunctions of the Hodge Laplacian $\Delta$ and set
$\mathcal{L}^{\Delta}_m=\textrm{span}\{ \xi_1, \ldots, \xi_m\}$
the vector space generated by these functions.
Clearly,   $\mathcal{H}_T^1(M)$ is a  subspace of $\mathcal{L}^{\Delta}_m$,
and by Lemma \ref{dim} $\dim \mathcal{H}_T^1(M) =2g+k-1$.
 
Next, we set an orthonormal basis  of  $C^{\infty}(M)$ 
given by eigenfunctions of the Jacobi operator, say $\{\phi_1,\phi_2,\ldots\}$.
We denote by  $\mathcal{J}^{p} :=\z \phi_1,\cdots,\phi_{p}\x^{\bot}$
the linear space orthogonal to the first $p$  eigenfunctions of the  Jacobi operator.

Initially, we look for vector fields $\xi\in \mathcal{L}^{\Delta}_m$  such that the functions 
$w_i,\bar w_i\in \mathcal{J}^{\alpha-1}$, for some $\alpha\in \mathbb N$ and 
$i\in\{1,2, 3\}.$ 
It is equivalent to find a solution to the following system with $6(\alpha-1)$ 
homogenous linear equations in the variable $\xi$
\begin{equation}\label{sys}
\int_Mw_i\phi_k\, dM=\int_M\bar w_i\phi_k\, dM=0,
\end{equation}
$1\leq i\leq 3$ and $1\leq k\leq \alpha-1$.
In particular, if $m(\alpha):=\dim \ \mathcal{L}^{\Delta}_m>6(\alpha-1),$ 
then the system  (\ref{sys}) has at least one  non trivial solution 
$\xi\in \mathcal{L}^{\Delta}_m$  such that $w_i,\bar w_i\in \mathcal{J}^{\alpha-1}$ 
for all $1\leq i\leq 3.$ By Courant minimax principle we have
\begin{eqnarray*}
\lambda_{\alpha}^J\int_M w_i^2\,dM \leq Q(w_i,w_i) 
\quad \mbox{and}\quad 
\lambda_{\alpha}^J\int_M\bar w_i^2 \,dM\leq Q(\bar w_i,\bar w_i). 
\end{eqnarray*}
Now, using Lemma \ref{lapwi} we get 
{\setlength\arraycolsep{1pt}
\begin{eqnarray*}
\lambda_{\alpha}^J\int_Mw_i^2\,dM&\leq &-4H^2\int_Mw_i^2\,dM+2H\int_M\z E_i,A\xi\x w_i\,dM\\
 &&+\int_M\z E_i,\Delta \xi\x w_i\,dM -2\int_Mg_i\z A,\nabla \xi\x w_i\,dM\\
 &&+\int_{\partial M}(w_i\eta(w_i)+II^{\partial W(N,N)}w_i^ 2)\,ds.
\end{eqnarray*}}

Summing up $i=1, 2, 3$ and using Lemma \ref{wi2} we obtain 
{\setlength\arraycolsep{1pt}
\begin{eqnarray*}
\lambda_{\alpha}^J\int_M\|\xi\|^2\, dM&\leq&-4H^2\int_M\|\xi\|^2\, dM+2H\int_M\z A\xi,\xi\x \, dM \\
&&+\int_M\z\Delta \xi,\xi\x\, dM+2\int_{\partial M}H^{\partial W}\|\xi\|^2 \, dM.
\end{eqnarray*}}

Aplying the same arguments to the test functions $\bar w_i$ we get
{\setlength\arraycolsep{1pt}
 \begin{eqnarray*}
\lambda_{\alpha}^J\int_M\|\xi\|^2\, dM&\leq &-4H^2\int_M\|\xi\|^2\, dM+2H\int_M\z A{\star\xi},\star\xi\x \, dM\\
&&+\int_M\z\Delta {\star\xi},\star\xi\x\, dM+2\int_{\partial M}H^{\partial W}\|\xi\|^2\,ds.
\end{eqnarray*}}

Then, summing these last two inequalities and noting that 
$\z A\xi,\xi\x+\z A{\star\xi},\star\xi\x=2H\|\xi\|^2$ we have 
{\setlength\arraycolsep{0pt}
 \begin{eqnarray}\label{lambda}
\lambda_{\alpha}^J\int_M\|\xi\|^2 dM&\leq
&+2\int_{\partial M}H^{\partial W}\|\xi\|^2 dM -2H^2\int_M\|\xi\|^2dM\\
&&+\frac{1}{2}\int_M(\z\Delta \xi,\xi\x+\z\Delta {\star\xi},\star\xi\x) dM. \nonumber
    \end{eqnarray}}
Finally,  if $\xi\in \mathcal{L}^{\Delta}_m$ we get  $\xi=\sum_i\alpha_i\xi_i$ and therefore
\begin{equation}\label{xi}
\int_M\z\Delta {\star\xi},\star\xi\x dM=\int_M\z\Delta \xi,\xi\x dM
=\lambda_{m(\alpha)}\int_M\|\xi\|^2 dM. 
    \end{equation}
Using (\ref{xi}) into (\ref{lambda}) and using the fact that $H^{\partial W}\leq 0$ we obtain
\[
\lambda_{\alpha}^J\leq -2H^2+\lambda_{m(\alpha)}^{\Delta},
\]
where $m(\alpha)>6(\alpha-1)$.  
It concludes the first part of Theorem \ref{esp}.

In order to get the lower bound for the index of $M$ 
we take $\{\phi_1,\phi_2,\cdots\}$  an orthonormal basis of the  space  $\mathcal{F}$ 
given by eigenfunctions of the Jacobi operator. 
From Lemma \ref{wi} we know that if $\xi\in\mathcal{H}_T^1(M)$, then the test 
functions  $w_1$, $w_2$ and $w_3$, belong to $\mathcal{F}$.

We look for  vector fields $\xi\in \mathcal{H}_T^1(M)$  such that for $1\leq i\leq3$,
the test functions  $w_i,\bar w_i\in \mathcal{J}^{\alpha-1}$, for some $\alpha\in \mathbb N,$ 
and $\bar w_i\in\mathcal{F}.$ 
In other words, we have the following system with $6\alpha-3$ homogeneous linear equations in 
the variable $\xi$:
\begin{equation}\label{sys2}
\int_M\bar w_i=\int_Mw_i\phi_k=\int_M\bar w_i\phi_k=0,
\end{equation}
where $1\leq i\leq 3$ and $1\leq k\leq \alpha-1.$

If $\dim \mathcal{H}_T^1(M)=2g+k-1>6\alpha-3,$ then the system (\ref{sys2})  has at least one 
non trivial solution $\xi\in \mathcal{H}_T^1(M)$.
Following the same steps as above we get 
{\setlength\arraycolsep{1pt}
 \begin{eqnarray*}
\lambda_{\alpha}^J\int_M\|\xi\|^2&\leq&-2H^2\int_M\|\xi\|^2.
    \end{eqnarray*}}

It implies that $\lambda_{\alpha}^J< 0$ and then $\indw (M)\geq \alpha.$ 
Since $\alpha$ can be chosen as the largest integer such that $2g+k-1>6\alpha-3$ we get 
\[
\indw(M)\geq \frac{2g+k-4}{6}.
\]


\bibliographystyle{amsplain}
\bibliography{bibliography}

\end{document}